\UseRawInputEncoding
\documentclass{amsart}
\usepackage{amsmath}
\usepackage{amssymb}
\newtheorem{thm}{Theorem}[section]
\newtheorem{cor}[thm]{Corollary}
\newtheorem{conj}[thm]{Conjecture}
\newtheorem{lem}[thm]{Lemma}
\newtheorem{prop}[thm]{Proposition}
\newtheorem{cons}[thm]{Construction}
\theoremstyle{definition}
\newtheorem{defn}[thm]{Definition}
\theoremstyle{remark}
\newtheorem{rem}[thm]{Remark}
\newtheorem{ex}[thm]{Example}
\newtheorem{exs}[thm]{Examples}
\long\def\Thm#1{\begin{thm} #1 \end{thm}}
\long\def\Cor#1{\begin{cor} #1 \end{cor}}

\long\def\Lem#1{\begin{lem} #1 \end{lem}}
\long\def\Prop#1{\begin{prop} #1 \end{prop}}

\def\Sect{\section}
\def\Rarr#1#2{\xrightarrow[#2]{#1}}

\long\def\Ref#1#2#3#4#5#6{
\bibitem{#1}
{\rm #2,}
\textit{#3.}
{\rm #4}
\textbf{#5}
{\rm #6.}
}
\long\def\Refb#1#2#3#4{
\bibitem{#1}
{\rm #2,}
\textit{#3.}
#4.
}
\def\Rr{{\mathbb R}}
\def\SS{{\mathcal S}}
\def\phi{\varphi}
\def\into{\hookrightarrow}
\def\iso{\cong}
\def\leq{\leqslant}
\def\geq{\geqslant}
\def\supp{{\rm supp}}

\def\st{\mid}
\def\Zero{{\rm Zero}}
\def\Hom{{\rm Hom}}

\begin{document}

\title{On covering dimension and sections of \goodbreak vector bundles}

\author{M.~C.~Crabb}
\address{%
Institute of Mathematics\\
University of Aberdeen \\
Aberdeen AB24 3UE \\
UK}

\email{m.crabb@abdn.ac.uk}
\date{February 2023, revised March 2024}
\begin{abstract}
An elementary result in point-set topology is used,
with knowledge of the mod $2$ cohomology of real projective spaces,
to establish
classical results of Lebesgue and Knaster-Kuratowski-Mazurkiewicz,
as well as the topological central point theorem of Karasev, which is
applied to deduce results of Helly-Lov\'asz, B\'ar\'any and Tverberg.
\end{abstract}
\maketitle
\Sect{Introduction}
Throughout this note,
$X$ will be a compact Hausdorff topological space
and $\xi_1$, $\ldots$, $\xi_n$ ($n\geq 1$) 
will be $n$ finite-dimensional real vector bundles over $X$.
We write $\xi =\xi_1\oplus\cdots\oplus\xi_n$.

\smallskip

Suppose that $A_1,\ldots , A_n$ are closed subspaces
covering $X$ and that, for each $k\in\{ 1,\ldots ,n\}$,
the restriction of $\xi_k$ to $A_k$ admits
a nowhere zero section, $s_k'$ say. 
By Tietze's theorem, $s_k'$ extends to a
section $s_k$ of $\xi_k$ on $X$. Then $s=(s_1,\ldots ,s_n)$
is a nowhere zero section of $\xi$.
(And conversely, if $\xi$ admits a nowhere zero section
$s=(s_k)$, we may construct such closed subspaces $A_k$
by choosing an inner product on each $\xi_k$ and
setting $A_k=\{ x\in X \st \| s_k(x)\| = \max\{ \| s_j(x)\| \st
1\leq j\leq n\}\}$.)
The main result of this paper, 
modelled on a cohomological lemma \cite[Lemma 3.2]{toric}
of Karasev and stated as Theorem \ref{palais_thm},
is in the same vein as this classical observation.  The methods
are from elementary point-set topology. 
As applications we derive in Sections 3 and 4,
using ideas introduced by Karasev in \cite{tcp, toric},
a classical result of Lebesgue
and Knaster-Kuratowski-Mazurkiewisz and the more recent topological
central point theorem of Karasev with, as corollaries, results of
Helly-Lov\'asz and B\'ar\'any \cite{barany}.

\smallskip

When we discuss Euler classes, we shall use representable cohomology 
as, for example, in \cite[Section 8]{borsuk}.

\smallskip

From the foregoing summary it should be clear that most of the ideas
presented in this note derive from the paper \cite{toric} of Karasev.
It is hoped, nevertheless, that the elementary approach taken here
may have some conceptual advantages.
\Sect{The principal result}
\Thm{\label{palais_thm}
Let $\xi_1$, $\ldots$, $\xi_n$ be $n$ finite-dimensional
real vector bundles over
a compact Hausdorff topological space $X$. 
Suppose that $(U_i)_{i\in I}$ is a finite open cover of $X$
such that each point of $X$ lies in $U_i$ for at most $n$ 
indices $i\in I$.

Suppose that for each $i\in I$ and $k\in\{ 1,\ldots ,n\}$ there
exists a section of $\xi_k$ with no zeros in $U_i$.
Then $\xi=\xi_1\oplus\cdots\oplus\xi_n$ 
admits a global nowhere zero section.
}

\begin{proof}
We begin with an argument from \cite[Lemma 2.4]{palais}. 
Choose a partition of unity $(\phi_i)_{i\in I}$ subordinate to the
cover. For $x\in X$, define
$$
J(x) =\{ j\in I \st \phi_j(x)=\max\{ \phi_i(x) \st i\in I\}\}.
$$
By assumption $\# J(x)\leq n$ (and $J(x)$ is non-empty).
For a non-empty subset $J\subseteq I$, we now set
$$
U_J=
\{ x\in X \st \text{\ for all
$j\in J$, $\phi_j(x)>0$ and $\phi_i(x)<\phi_j(x)$ for all $i\in I-J$}\}.
$$
It is clearly an open subset of $\bigcap_{j\in J}U_j$.
Moreover, $x\in U_{J(x)}$, so that the sets $U_J$ cover $X$.
If $J$ and $J'$ are distinct subsets of $I$ with $\# J=\# J'$,
then $U_J\cap U_{J'}=\emptyset$. (For there exist elements $j\in J$,
with $j\notin J'$, and $j'\in J'$, with $j'\notin J$.)

Choose a partition of unity $(\psi_J)$ subordinate to the open
cover $(U_J)_{J\subseteq I\, :\, 1\leq \# J\leq n}$ of $X$ and, 
for each $J$ with $\# J=k\leq n$,
a section $s_J$ of $\xi_k$ which is nowhere zero on $U_J$
(possible because $U_J\subseteq U_i$ if $i\in J$).
We can then define a section $s_k$ of $\xi_k$ by
$$
s_k(x) = \sum_{\# J=k} \psi_J (x) s_J(x).
$$
Notice that, if $\psi_J(x)\not=0$ for some $J$ with
$\# J=k$, then $s_k(x)\not=0$. For, if $J'$ is a different
subset with $\# J'=k$, $\psi_{J'}(x)=0$.

The section $s=(s_1,\ldots ,s_n)$ of $\xi$ is nowhere zero.
\end{proof}
\Cor{\label{closed_cover}
Suppose that $(A_i)_{i\in I}$ is a finite closed cover of $X$
such that each point of $X$ lies in $A_i$ for at most $n$ 
indices $i\in I$.

Suppose that for each $i\in I$ and $k\in\{ 1,\ldots ,n\}$ the
restriction of $\xi_k$ to $A_i$ admits a nowhere zero section.
Then $\xi$ admits a global nowhere zero section.
}
\begin{proof}
It is an elementary exercise to
show that there is an open cover $(U_i)_{i\in I}$
such that (i) $A_i\subseteq U_i$, (ii) for each $(i,k)$
there exits a global section of $\xi_k$ with no
zeros in $U_i$, and (iii) each point of $X$ lies in at most $n$
of the sets $U_i$. Then we can apply Theorem \ref{palais_thm}.

(Here are the details when $\# I >n$.
Consider a subset $J\subseteq I$ with $\# J>n$.
The open sets $X-A_j$, $j\in J$, cover $X$. Choose a partition
of unity $(\chi_j)_{j\in J}$ subordinate to this cover.
Then we can define $U_j^J =\{ x\in X \st \# J\cdot \chi_j(x)< 1\}$. 
By construction $\bigcap_{j\in J} U_j^J=\emptyset$ (because
$\sum_{j\in J}\chi_j(x)=1$ for any $x$)
and $A_j\subseteq U_j^J$ for $j\in J$.

Every nowhere zero section of $\xi_k$ over the closed set $A_i$
extends to a global section of $\xi_k$ and such a section will be
nowhere zero on an open neighbourhood of $A_i$. So it is easy to
choose $U_i$ in the intersection of the sets $U_i^J$ with
$\# J>n$ and $i\in J$ to satisfy (i) and (ii).
And then, for any set $J$ with $\# J >n$, we have
$\bigcap_{i\in J}U_i =\emptyset$.)
\end{proof}
\Sect{A theorem of Lebesgue}
If the mod $2$ cohomology Euler class $e(\xi )$ of $\xi$ is non-zero,
the {\it zero-set} $\Zero (s)=\{ x\in X\st s(x)=0\}$ 
of any section $s$ of $\xi$ is non-empty.
\Prop{\label{lift}
Let $p: X\to Y$ be a continuous map from 
$X$ to a compact Hausdorff space  $Y$. Suppose that $(B_i)_{i\in I}$
is a finite cover of $Y$ by closed sets such that any
point of $Y$ lies in at most $n$ of the sets $B_i$.

If the mod $2$ cohomology Euler class $e(\xi )$ of $\xi$
is non-zero, then there exist $i\in I$ and $k\in\{ 1,\ldots ,n\}$
such that $p(\Zero (s))\cap B_i$ is non-empty
for each section $s$ of $\xi_k$.
}
\begin{proof}
We take $A_i =p^{-1}(B_i)$ in Corollary \ref{closed_cover}.
Since $e(\xi )\not=0$, $\xi$ does not admit a nowhere zero
section. Hence there is a pair $(i,k)$ such that
$\Zero (s)\cap B_i\not=\emptyset$ for every section
$s$ of $\xi_k$.
\end{proof}
For a finite set $V$ of cardinality $n+1$, we write
$\Rr [V]$ for the $(n+1)$-dimensional real vector space
of maps $t: V \to \Rr$, $\Delta (V)$ for the
$n$-simplex of maps $t$ such that $t(v)\geq 0$ for all $v\in V$ and
$\sum_v t(v)=1$,
$S(\Rr [V])$ for the unit $n$-sphere in $\Rr [V]$
of vectors $t$ with $\sum_v t(v)^2=1$ and
$P(\Rr [V])$ for the $n$-dimensional real projective space of lines
$[t]$ in $\Rr [V]$ (where $t\in\Rr [V]$ is non-zero).
There is a surjective map
$$
\pi_V : P(\Rr [V]) \to \Delta (V)
$$
defined by $(\pi_V [t])(v) =t(v)^2$ for $t\in S(\Rr [V])$.
There is also an embedding
$$
\sigma_V : \Delta (V) \into S(\Rr [V])\, ,
$$
defined by $(\sigma_V (t)(v) = \sqrt{t(v)}$,
such that the composition
$$
\Delta (V) \Rarr{\sigma_V}{} S(\Rr [V]) \Rarr{}{}
P(\Rr [V]) \Rarr{\pi_V}{} \Delta (V)
$$
is the identity.
The Hopf line bundle over $P(\Rr [V])$, with fibre
$\Rr t\subseteq \Rr [V]$ at $[t]$, is denoted by $H$.

The {\it support}, $\supp (t)$, of $t\in\Rr[V]$
is the set of points $v\in V$ such that $t(v)\not=0$.
For an integer $d$, $0\leq d\leq \# V$, we write
$\SS_d(V)$ for the set of finite subsets $T$
of $\Delta (V)$ such that any two elements of $T$
have disjoint supports and $\# V -\# T=d$.
For $T\in\SS_d(V)$, we write $\Delta_T$ for the
convex hull of $T$; it is a simplex of codimension
$d$ in $\Delta (V)$.
\Lem{\label{sec}
For $T\in\SS_d(V)$, the simplex $\Delta_T$ can be
expressed as $\pi_V(\Zero (s))$ for some section
$s$ of the vector bundle $dH=\Rr^d \otimes H$.
}
\begin{proof}
Let $E_T$ be the codimension $d$ vector subspace of $\Rr [V]$
spanned by $\sigma_V(T)$.
The section $s_T$ of $\Hom (H,\Rr [V]/E_T)$ over
$P(\Rr [V])$  given by the projection
$s_T([t])  = \{ [t]\into \Rr [V]\to \Rr [V]/E_T\}$,
for $t\in S(\Rr [V])$, has
the property that $\pi_V(\Zero (s_T))=\Delta_T$.
Choose some isomorphism $\Rr [V]/E_T\iso \Rr^d$ to
get the required section $s$.
\end{proof}
\Thm{\label{lkkm}
Let $(V_l)_{l=1}^m$ be a family of $m$ finite sets with
$\# V_l = d_ln_l+1$, where $d_l\geq 1$ and $n_l\geq 1$,
$l=1,\ldots ,m$,
are positive integers. Write $n=n_1+\ldots +n_m$.
Suppose that $(B_i)_{i\in I}$ is a finite closed 
cover of
$$
Y= \Delta (V_1)\times\cdots\times \Delta (V_m)
$$
such that any point of $Y$ lies in at most $n$ of
the sets $B_i$.

Then for some $i\in I$ and $l\in\{ 1,\ldots ,m\}$ the projection
of $B_i$ to the $l$th factor $\Delta (V_l)$ meets each of the
codimension $d_l$ simplices $\Delta_T$ for $T\in\SS_{d_l}(V_l)$.
}
The Lebesgue theorem \cite[Theorem 4.1]{toric} is the special case 
$n_l=1$, $d_l=1$, for all $l$;
the case $m=1$, $d_1=1$ is a result of 
Knaster, Kuratowski, Mazurkiewicz \cite[Remark 2.2]{toric}.

We follow the method of Karasev \cite[Theorems 2.1 and 4.1]{toric}.
\begin{proof}
This is an immediate consequence of Proposition \ref{lift}.
Take $p$ to be the product
$$
\pi_{V_1} \times\ldots\times\pi_{V_m} :
X=P(\Rr [V_1]) \times \cdots \times P(\Rr [V_m])
\to Y=\Delta (V_1)\times\cdots\times\Delta (V_m)
$$
and $\xi_k$ to be the multiple $d_lH_l$ 
of the Hopf bundle $H_l$ over
$P(\Rr [V_l])$ if $n_1+ \ldots +n_{l-1}<k\leq
n_1+\ldots +n_l$.
The Euler class $e(\xi ) =e(H_1)^{d_1n_1}\ldots e(H_m)^{d_mn_m}$
is non-zero.
\end{proof}
The next lemma gives a way of checking that $e(\xi )$ in
the application of Proposition \ref{lift} is non-zero.
\Lem{\label{zero}
Suppose that $X$ is a closed subspace of a compact Hausdorff space
$\hat X$ and that $\xi$ is the restriction of a vector bundle 
$\hat\xi$ over $\hat X$.
Let $\eta$ be a real vector bundle over $\hat X$ with the property
that the restriction of $\eta$ to each connected component 
of the complement $\hat X-X$ admits a nowhere zero section.

If the mod $2$ cohomology 
Euler class $e(\hat \xi \oplus \eta)$ is non-zero,
then $e(\xi)$ is non-zero.
}
\begin{proof}
Fix a Euclidean metric on $\eta$.
For each component $C$ of $\hat X-X$ choose a continuous
section $s_C$ of the sphere bundle
$S(\eta\,|\, C)$, and choose a continuous function $\rho : \hat X\to
[0,1]$ such that $\rho^{-1}(0)=X$.
Then one can define a continuous section $s$ of $\eta$ 
with zero-set $\Zero (s)=X$ by
$s(x)=\rho (x)s_C(x)$ if $x\in C$, 
$s(x)=0$ if $x\in X$.

Since $e(\hat\xi\oplus\eta )=e(\hat \xi )\cdot e(\eta)$ is non-zero,
the restriction $e(\xi)$ of $e(\hat \xi )$ to $\Zero (s)=X$
is non-zero. See, for example, \cite[Proposition 2.7]{borsuk}.
\end{proof}
This allows us to deduce Karasev's strengthened KKM theorem
\cite[Theorem 2.1]{toric}).
\Prop{
Let $(B_i)_{i\in I}$ be closed subsets of a simplex $\Delta (V)$ 
with vertex set $V$ of cardinality $dn+r+1$, 
where $d\geq 1$ and $r\geq 0$ are integers,
such that any point of $Y=\bigcup_{i\in I} B_i$ lies in at most $n$ 
of the sets $B_i$.
Then,
either for some $i\in I$ the subset $B_i$ meets each codimension
$d$ simplex $\Delta_T$ for $T\in \SS_d(V)$,
or
some connected component of $\Delta (V)-Y$ intersects
every $dn$-dimensional simplex $\Delta_T$ for $T\in \SS_r(V)$.
}
\begin{proof}
Considering $\pi_V :\hat X=P(\Rr [V])\to \hat Y=\Delta (V)$, take
$X\subseteq \hat X$ to be $\pi_V^{-1}(Y)$ and $p:X\to Y$ to be
the restriction of $\pi_V$.
Let each $\xi_k$ be the restriction of $dH$,
so that $\xi$ is the restriction of $\hat\xi = dnH$
to $X$. Take $\eta =rH$,
so that $e(\hat\xi\oplus \eta )=e(H)^{dn+r}$ is
non-zero. 

If, for each component of $\Delta (V)-Y$,
written as the image $\pi_V(C)$ of a component $C$
of $\hat X-X$,
there is some simplex $\Delta_T$, where $T\in\SS_r(V)$, 
such that $\pi_V(C)\cap\Delta_T=\emptyset$,
then Lemma \ref{sec} provides a nowhere zero section of $rH=\eta$
over $C$.
By Lemma \ref{zero}, $e(\xi )$ is then non-zero, and
we can apply Proposition \ref{lift} to deduce the
existence of an $i\in I$ such that $B_i$ meets
each codimension $d$ simplex.
\end{proof}
\Sect{Karasev's topological central point theorem}
An early result of the following type appears in 
\cite[Lemma 3.1]{yang}.
\Prop{\label{yang_prop}
Let $f : X\to 
Z$ be a continuous map from $X$ to
a compact Hausdorff space $Z$ with covering dimension less than $n$.
Suppose that the mod $2$ cohomology class $e(\xi )$
is non-zero.

Then there exists a point $z\in Z$ and $k\in\{ 1,\ldots ,n\}$
such that $z\in f(\Zero (s))$ for each section $s$ of $\xi_k$.
}
\begin{proof}
Suppose that for each point $z\in Z$ there exist sections $s_k^z$
of $\xi_k$, $1\leq k\leq n$, such that 
$z\notin f(\Zero (s_k^z))$ for each $k$. Then the open sets
$(Z-f(\Zero (s_1^z))\cap \cdots \cap (Z-f(\Zero (s_n^z))$,
$z\in Z$, cover $Z$.
Since $Z$ is compact with covering dimension $<n$, this
open cover may be refined by a finite open cover
$(W_i)_{i\in I}$ such that each point of $Z$ lies in at most
$n$ of the sets $W_i$.

Set $U_i=f^{-1}(W_i)$. Then we may apply Theorem \ref{palais_thm}
to the open cover $(U_i)_{i\in I}$ of $X$ to conclude that
there exist $i$ and $k$ such that for every section $s$ of
$\xi_k$ the zero-set $\Zero (s)$ meets $U_i$. 
So $W_i\not\subseteq Z-f(\Zero (s))$ for every section
$s$ of $\xi_k$. But $W_i\subseteq Z-f(\Zero (s_k^z))$ for
some $z\in Z$. This contradiction completes the proof.
\end{proof}
\Thm{\label{centrept}
Let $(V_l)_{l=1}^m$ be a family of $m$ finite sets with
$\# V_l = d_ln_l+1$, where $d_l,\, n_l\geq 1$,
$l=1,\ldots, m$ are positive integers. Write $n=n_1+\ldots +n_m$.
Suppose that
$$
g: Y=\Delta (V_1)\times\cdots\times \Delta (V_m)\to Z
$$
is a continuous map to a compact Hausdorff space $Z$ with covering
dimension less than $n$. 

Then for some $l\in\{ 1,\ldots ,m\}$
$$
\bigcap_{T\in\SS_{d_l}( V_l )} g(\prod_{j<l} \Delta (V_j)\times\Delta_T
\times\prod_{j>l} \Delta (V_j)) \not=\emptyset\, .
$$
}
Karasev's topological central point theorem, as in
\cite[Theorem 1.1]{tcp} and 
\cite[Theorem 5.1]{toric}, is the special case $m=1$.
\begin{proof}
We apply Proposition \ref{yang_prop} with $X$ and $\xi_k$
as in Theorem \ref{lkkm} and with $f$ equal to the composition of
$$
\pi_{V_1}\times\cdots\times\pi_{V_m}:
X=P(\Rr [V_1])\times\cdots\times P(\Rr [V_m])\to
Y=\Delta (V_1)\times\cdots\times \Delta(V_m)
$$
with $g: Y\to Z$.
We recall that
$\xi_k= d_lH_l$ if $n_1+\ldots +n_{l-1}<k\leq n_1+\ldots +n_l$,
so that the Euler class 
$e(\xi )= e(H_1)^{d_1n_1}\cdots e(H_m)^{d_mn_m}$ is non-zero.
\end{proof}
As an application we prove a result of Helly-Lov\'asz 
\cite[Theorem 3.1]{barany}.
\Cor{Suppose that $C_{l,v}$, $l=1,\ldots ,m$,
$v\in V_l$, $\# V_l=d_l+1$, are convex subsets of a real vector space
$E$ with the property that the intersection
$C_{1,v_1}\cap \cdots \cap C_{m,v_m}$ is non-empty for each
$(v_1,\ldots ,v_m)\in V_1\times\cdots\times V_m$.

If $\dim E <m$,
then, for some $l$, the intersection $\bigcap_{v\in V_l} C_{l,v}$
is non-empty.
}
\begin{proof}
For each $(v_1,\ldots ,v_m)\in V_1\times\cdots\times V_m$ choose
$z_{v_1,\ldots ,v_m}\in C_{1,v_1}\cap \cdots C_{m,v_m}$.
We apply Theorem \ref{centrept} with $n=m$, $n_l=1$, and 
$Z\subseteq E$ the convex hull of the points $z_{v_1,\ldots ,v_m}$.
Take $g$ to be the piecewise linear map
$$
(\sum_{v_1\in V_1} t_1(v_1),\ldots ,\sum_{v_m\in V_m} t_m(v_m))
\mapsto \!\!\!\sum_{(v_1,\ldots ,v_m)\in V_1\times\cdots\times V_m}
\!\!
t_1(v_1)\cdots t_m(v_m)\, z_{v_1,\ldots ,v_m}\in Z.
$$
We conclude from Theorem \ref{centrept},
noting that a codimension $d_l$ simplex in $\Delta (V_l)$
is a point,
that there is some $l$ and $z\in Z$
such that, for each $v\in V_l$ the vector $z$ can be written as
$$
z=\sum_{(v_1,\ldots ,v_m)\in V_1\times\cdots\times V_m\, :\, v_l=v}
t_1(v_1)\cdots t_m(v_m)\, z_{v_1,\ldots ,v_m},
$$
where $\sum_{v_j\in V_j} t_j(v_j)=1$ for each $j$
and so $\sum_{(v_1,\ldots ,v_m)} t_1(v_1)\cdots t_m(v_m)=1$.
Since each $z_{v_1,\ldots ,v_m}$ with $v_l=v$ lies in the
convex set $C_{l,v}$, we see that $z\in C_{l,v}$, as
required.
\end{proof}
B\'ar\'any's dual result \cite[Theorem 2.1]{barany}
(as formulated in \cite[Theorem 3.1]{HK} and \cite[Theorem 3]{MR})
can be obtained in a similar fashion.
\Cor{\label{barany_cor}
Let $K\subseteq E$ be a non-empty compact convex subspace of 
a finite-dimensional
real vector space $E$.
Suppose that $V_1,\ldots ,V_m$ are finite sets
with $\# V_l=d_l+1$, $l=1,\ldots ,m$, and that $\phi_l :V_l \to E$
are maps with the property that for
each $(v_1,\ldots ,v_m)\in V_1\times\cdots\times V_m$ the
convex hull of $\{ \phi_1(v_1),\ldots ,\phi_m(v_m)\}$ in $E$ 
is disjoint from $K$.

If $\dim E<m$, then, for some $l\in\{ 1,\ldots ,m\}$, 
the convex hull of $\phi_l(V_l)$
in $E$ is disjoint from $K$.
}
\begin{proof}
Choose a basepoint $*\in K$.
Let $A$ to be the affine space of affine linear maps
$z: E\to\Rr$ such that $z(*)=-1$. Notice that the dimension
of $A$, as affine space, is equal to $\dim E$.

For each $(v_1,\ldots ,v_m)\in V_1\times\cdots\times V_m$ choose
an affine linear map $z_{v_1,\ldots ,v_m}\in A$ taking strictly positive values on $\{ \phi_1(v_1),\ldots ,\phi_m(v_m)\}$
and strictly negative values on $K$.
We again apply Theorem \ref{centrept} with $n=m$, $n_l=1$, and
$Z\subseteq A$ the convex hull of the points $z_{v_1,\ldots ,v_m}$.
As in the proof of Helly's theorem,
take $g$ to be the piecewise linear map
$$
(\sum_{v_1\in V_1} t_1(v_1),\ldots ,\sum_{v_m\in V_m} t_m(v_m))
\mapsto \!\!\!\sum_{(v_1,\ldots ,v_m)\in V_1\times\cdots\times V_m}
\!\!
t_1(v_1)\cdots t_m(v_m)\, z_{v_1,\ldots ,v_m}\in Z.
$$
Theorem \ref{centrept} provides some $l$ and $z\in Z$
such that, for each $v\in V_l$ the affine linear map
$z$ can be written as
$$
z=\sum_{(v_1,\ldots ,v_m)\in V_1\times\cdots\times V_m\, :\, v_l=v}
t_1(v_1)\cdots t_m(v_m)\, z_{v_1,\ldots ,v_m},
$$
where $\sum_{(v_1,\ldots ,v_m)} t_1(v_1)\cdots t_m(v_m)=1$,
so that $z$ takes a strictly positive value at each 
$\phi_l(v)\in \phi_l(V_l)$
and strictly negative values on $K$. 
\end{proof}
As observed by Sarkaria \cite{sarkaria} (and expounded in
\cite{soberon}), Tverberg's theorem is an easy consequence
of Corollary \ref{barany_cor}.
The following generalization,
discussed in \cite[Theorem 3.8]{soberon} and
due to Arocha, B\'ar\'any, Bracho, Fabila and Montejano \cite{ABBFM},
can be viewed as a coincidence theorem.
\Cor{
Let $r\geq 0$ and $m\geq 1$ be integers. 
For $l=1,\ldots ,m$, $s=0,\ldots ,r$, let $V_{l,s}$
be non-empty finite sets
and $\phi_{l,s} : V_{l,s}\to F$ be maps to a finite-dimensional
real vector space $F$ satisfying the two conditions:

\smallskip

\par\noindent
{\rm (i)}
for each $l\in\{1,\ldots ,m\}$, there is a non-zero vector in $F$
that can be expressed, for each $s=0,\ldots ,r$,
as a linear combination with non-negative
coefficients of the elements of $\phi_{l,s}(V_{l,s})$;
\par\noindent
{\rm (ii)}
for each $s$ and $v_l\in V_{l,s}$, $l=1,\ldots ,m$, the convex hull
of $\{ \phi_{1,s}(v_1),\ldots ,\phi_{m,s}(v_m)\}$
is disjoint from $\{ 0\}$.

\smallskip

Then, if $r\cdot\dim F <m$,
there is a partition $\{ 1,\ldots ,m\}=\bigsqcup_{s=0}^r I_s$
into $r+1$ non-empty subsets $I_s$
and a non-zero vector $c\in F$ 
such that
$$
\sum_{i\in I_s} \lambda_i\phi_{i,s}(v_i) =c
\text{\quad for $s=0,\ldots ,r$,}
$$
for some $\lambda_i\geq 0$
and $v_i\in V_{i,s}$ for $i\in I_s$.

If, further, there is some 
affine hyperplane $H$ in $F$ that contains
all the subsets $\phi_{l,s}(V_{l,s})$ but
does not contain $0$,
then $c$ may be chosen in the hyperplane $H$
and then $\sum_{i\in I_s}\lambda_i=1$ for each $s$.
}
\begin{proof}
Let $L_r$ be the quotient of $\Rr^{r+1}=\bigoplus_{s=0}^r \Rr e_s$ by
the subspace generated by $e_0+\ldots +e_r$ and write $[e_s]$
for the coset of $e_s$.

We take $V_l=\bigsqcup_{s=0}^r V_{l,s}$, $E=L_r\otimes F$,
$\phi_l:V_l \to E$ defined by $\phi_l(v)=[e_s]\otimes \phi_{l,s}(v)$
for $v\in V_{l,s}$, and apply Corollary \ref{barany_cor} with
$K=\{ 0\}$.
By assumption, the convex hull of each $\phi_l(V_l)$ in $E$
contains $0$. (Notice that, if $a_1,\ldots ,a_r\in F$, then
$\sum_s [e_s]\otimes a_s=0$ if and only if $a_1=a_2=\cdots =a_r$.)

So there exist $v_i\in V_i$, $i=1,\ldots ,m$, $\lambda_i\geq 0$,
and a non-zero $c\in F$ such that
$\sum_{i\, :\, v_i\in V_{i,s}} \lambda_i\phi_{i,s}(v_i) =c$.
Take $I_s =\{ i\st v_i\in V_{i,s}\}$.

If there is a linear form $\alpha : F\to\Rr$ taking the value $1$ on
all $\phi_{l,s}(V_{l,s})$, we can scale 
to arrange that $\alpha (c)=1$,
and then $\sum_{i\in I_s}\lambda_i=1$.
\end{proof}
The original Tverberg theorem is the case in which $V_{l,s}=\{ *\}$
is a single point
for all $l,\, s$ and $\phi_{l,s}(*)$ is independent of $s$.
\begin{appendix}
\Sect{Cohomology}
It is a classical result that,
if $(A_k)_{k=1}^n$ is a closed cover of a compact Hausdorff space
$X$ and, for each $k$, $e_k$ is a mod $2$ cohomology class of 
$X$ that restricts to zero on $A_k$, then the product
$e_1\cdots e_n$ is zero.
Here is the corresponding version of Corollary \ref{closed_cover},
which was used by Karasev in the form \cite[Lemma 3.2]{toric}.
\Thm{\label{cohomology}
Let $X$ be a compact Hausdorff space and let $e_1,\ldots ,e_n$
by classes in the mod $2$ cohomology of $X$.

Suppose that $(A_i)_{i\in I}$ is a finite closed cover of $X$
such that each point of $X$ lies in $A_i$ for at most $n$ 
indices $i\in I$ and
that for each $i\in I$ and $k\in\{ 1,\ldots ,n\}$ the
restriction of $e_k$ to $A_i$ is zero.

Then $e_1\cdots e_n=0$.
}
\begin{proof}
By the argument used in the proof of Corollary \ref{closed_cover}
one can manufacture an open cover $(U_i)$ such that each
cohomology class is represented by a map which is null 
(not just null-homotopic) on
$U_i$. The construction in the proof of
Theorem \ref{palais_thm} gives an open cover $(U_J)$ indexed by
the non-empty subsets $J$ of $I$ with $\# J\leq n$.
The map representing $e_k$ is null on the disjoint union of
the sets $U_J$ with $\# J=k$. Since these $n$ open sets
cover $X$, the product $e_1\cdots e_n$
is represented by the null map, and so 
the cohomology class $e_1\cdots e_n$ is zero.
\end{proof}
\Sect{Fibrewise joins}
The principal result, Theorem \ref{palais_thm}, extends from sphere 
bundles to fibre bundles (understood to be locally trivial).
\Lem{\label{global}
Let $E\to X$ be a fibre bundle over a compact
Hausdorff space $X$ with each fibre a compact ENR. Then there is a
fibrewise embedding $j:E\into \Omega$ into an open subspace of
a trivial real vector bundle $X\times V$ admitting a fibrewise
retraction $r: \Omega\to E$.
}
\begin{proof}
We recall the well known argument (from, for example,
\cite[II, Lemma 5.8]{FHT}).
Choose a finite open cover $(U_i)_{i=1}^n$ of $X$,
with a partition of unity $(\phi_i)$ subordinate to the cover,
and trivializations of the restriction of $E$ to each subspace $U_i$:
$E\, |\, U_i \to U_i\times F_i: y\mapsto (x,f_i(y))$,
for a point $y\in E_x$ in the fibre of $E$ at $x\in X$,
where $F_i$ is a compact Euclidean Neighbourhood Retract (ENR)
embedded as a subspace 
$F_i\subseteq \Omega_i\subseteq V_i$
of an open subspace $\Omega_i$ of a Euclidean space $V_i$
with a retraction $r_i : \Omega_i\to F_i$.

Putting $V=\bigoplus_{i=1}^n V_i$,
define an embedding $j: E\to X\times V$ by 
$$
j(y)=(x,\phi_1(x)f_1(y),\ldots ,\phi_n(x)f_n(y))
$$ 
for $y\in E_x$.

Let $W_i\subseteq X\times V$ be the open subset of points
$(x,(v_j))$ such that either $n\phi_i(x)<\frac{1}{2}$ or
$$
\phi_i(x)>0\text{\ and\ }v_i/\phi_i(x)\in\Omega_i.
$$
Observe that $j(E)\subseteq W_i$.

Define $q_i :W_i\to X\times V$ by
$q_i(x,(v_j))=(x,(v_j))$ if $n\phi_i(x)<\frac{1}{2}$,
and
$$\textstyle
q_i(x,(v_j))=(x,(t\phi_j(x)f_j(y)+(1-t)v_j))
\text{\  if $n\phi_i(x)\geq \frac{1}{2}$,} 
$$
where $r_i(v_i/\phi_i(x))=f_i(y)$ and
$t=\min \{1,\, n\phi_i(x)-1/2\}$.
So, if $n\phi_i(x)\geq 1$, we have $t=1$
and $q_i(x,(v_j))\in j(E)$.
And, because $\sum_j\phi_j(x)=1$, there is at least
one $i$ such that $n\phi_i(x)\geq 1$.

Now take $\Omega$ and $r$ to be the open subset
$$
\Omega =\{ (x,(v_j))\in W_1\st q_i (q_{i-1} (\cdots
q_1(x,(v_j))\cdots ))\in W_{i+1},\, i=1,\ldots ,n-1\}
$$
and retraction
$r(x,(v_j))=q_n(q_{n-1}( \cdots q_1(x,(v_j))\cdots ))$.
\end{proof}
\Thm{\label{palais_thm_gen}
Let $E_1$, $\ldots$, $E_n$ be $n$ fibre bundles, with each fibre
a compact ENR, over
a compact Hausdorff topological space $X$. 
Suppose that $(A_i)_{i\in I}$ is a finite closed cover of $X$
such that each point of $X$ lies in $A_i$ for at most $n$ 
indices $i\in I$.

Suppose that for each $i\in I$ and $k\in\{ 1,\ldots ,n\}$ there
exists a section of $E_k \, | \, A_i$ over $A_i$.
Then the fibrewise join  $E=E_1*_X \cdots *_X E_n$ 
admits a global section.
}
We shall write points of the join $F_1*\cdots *F_n$ of
spaces $F_k$, $k=1,\ldots ,n$,
as $[(y_1,\ldots ,y_n),(t_1,\ldots ,t_n)]$,
where $y_k\in F_k$, $t_k\in [0,1]$, and $\sum t_k=1$.
\begin{proof}
The proof of Corollary \ref{closed_cover} using Theorem 
\ref{palais_thm} is readily adapted.

First, using Lemma \ref{global} to see that
a section of $E_k$ over $A_i$ extends to a section
over an open neighbourhood
of $A_i$ in $X$, we construct an open cover $(U_i)$ such that
(i) $A_i\subseteq U_i$, (ii) for each $(i,k)$ there is a section
of $E_k\, |\, U_i$ over $U_i$, and (iii) each point of $X$
lies in at most $n$ of the sets $U_i$.

Then having produced the open subsets $U_J$ as in the proof of
Theorem \ref{palais_thm} and chosen $\psi_J$, we choose 
sections $s_J$ of $E_k\, |\, U_J$, where $\# J=k$.
These combine over the finite disjoint union of the subsets
$U_J$ to give a section $s_k$ of $E_k$ over
$\bigsqcup_{\# J=k} U_J$.
A section $s$ of the fibrewise join $E_1 *_X\cdots *_X E_n$
is given by $s(x)=[(s_1(x),\ldots ,s_n(x)),(t_1(x),\ldots ,t_n(x))]$,
where $t_k(x)=\sum_{\# J=k} \psi_J(x)$, so that 
$t_1(x)+\cdots +t_n(x)=1$.
\end{proof}

\end{appendix}

\end{document}